# STATIONARY PROCESSES WHOSE FILTRATIONS ARE STANDARD

By X. Bressaud,[1] A. Maass,[2] S. Martinez[2] and J. San Martin[2]

*Université de la Mediterranée, Universidad de Chile, Universidad de Chile and Universidad de Chile*

We study the standard property of the natural filtration associated to a 0–1 valued stationary process. In our main result we show that if the process has summable memory decay, then the associated filtration is standard. We prove it by coupling techniques. For a process whose associated filtration is standard, we construct a product type filtration extending it, based upon the usual couplings and the Vershik's criterion for standardness.

**1. Introduction and notation.** Let $(X_n : n \leq 0)$ be a $\{0,1\}$-valued stationary process and $\mathcal{F}^X = (\mathcal{F}_n^X : n \leq 0)$ be its natural filtration, so $\mathcal{F}_n^X = \sigma(X_m; m \leq n)$.

DEFINITION 1. A filtration $\mathcal{F}$ is *standard* if it can be immersed on a filtration of diffusive product type (see [6, 7, 8, 15, 16]).

A necessary condition for $\mathcal{F}$ to be standard is that its tail $\mathcal{F}_{-\infty} = \bigcap_{n \leq 0} \mathcal{F}_n$ is trivial. But, as is shown by a counterexample in [15, 16], this condition is not sufficient.

In our main result we show that if $(X_n : n \leq 0)$ has (a slightly weaker condition than) summable memory decay, then $\mathcal{F}^X$ is standard. This is done in Theorem 3 of Section 3. For the proof, we construct explicitly a filtration $\mathcal{G} = (\mathcal{G}_n : n \leq 0)$, where $\mathcal{F}^X$ is immersed, and further, we show it is of diffusive product type. That is, there exists a sequence of i.i.d. uniform r.v.'s $(W_n : n \leq 0)$ such that $\mathcal{G} = \mathcal{F}^W$.

Received November 2004; revised August 2005.
[1]Supported by CNRS.
[2]Supported by Nucleus Millennium Information and Randomness P01-005 and P04-69F.

To be more precise, let $\Sigma = \{0,1\}^{-\mathbb{N}}$ be endowed with the law of $(X_n : n \leq 0)$. Let $(V_n : n \leq 0)$ be a sequence of i.i.d. r.v.'s uniformly distributed on $[0,1]$, independent of $\mathcal{F}^X$. We endow $[0,1]^{-\mathbb{N}}$ with the law of $(V_n : n \leq 0)$ and we fix the probability space $(\Omega, \mathcal{A}, \mathbb{P})$ as the product of above spaces, so $\mathbb{P}$ is the product of the laws of $(X_n : n \leq 0)$ and $(V_n : n \leq 0)$. On the other hand, the filtration $\mathcal{G} = (\mathcal{G}_n : n \leq 0)$ is given by $\mathcal{G}_n = \sigma(X_m, V_m : m \leq n)$. Clearly, $\mathcal{F}^X$ is immersed in $\mathcal{G}$ (see [6]). The above mentioned sequence $(W_n : n \leq 0)$ is constructed in Section 2.

The class of processes with summable memory decay has been studied in relation with regenerative representations and perfect simulation algorithms, in particular, see [2, 3, 5, 9]. Gibbs measures with Hölder potentials on fullshifts are examples of measures with summable memory decay (see [1, 13]); a rich discussion and a detailed list of relevant references on this class of measures can be found in [3, 9].

In Section 4 we assume $\mathcal{F}^X$ is standard and we construct an explicit diffusive product type extension $\mathcal{F}^U$ of $\mathcal{F}^X$.

**2. An independent sequence.** Let $n \leq 0$. We define $f_n = \mathbb{P}(X_n = 0 | \mathcal{F}_{n-1}^X)$ and

$$(1) \qquad W_n = f_n V_n \mathbf{1}(X_n = 0) + (1 - (1 - f_n)V_n)\mathbf{1}(X_n = 1),$$

where $\mathbf{1}(X_n = i)$ denotes the characteristic function of the event $\{X_n = i\}$, for $i = 0, 1$.

LEMMA 2. $(W_n : n \leq 0)$ is a sequence of i.i.d. r.v.'s uniformly distributed in $[0,1]$. Moreover, for all $n \leq 0$, $W_n$ is independent of $\mathcal{G}_{n-1}$, $\mathcal{G}_{n-1} \vee \sigma(W_n) = \mathcal{G}_n$, and $\mathcal{F}_{n-1}^X \vee \sigma(W_n) = \mathcal{F}_n^X \vee \sigma(V_n)$.

PROOF. First recall the following relation. Let $f$, $V$ and $Z$ be real bounded measurable functions and $\mathcal{B}$ be a sub $\sigma$-field such that $f$ is $\mathcal{B}$-measurable and $V$ is independent of $\mathcal{B} \vee \sigma(Z)$. Then, for any Borel real bounded function $h$, it holds $\mathbb{E}(h(fV)Z|\mathcal{B})(\omega) = \mathbb{E}(Z|\mathcal{B})(\omega) \int h(f(\omega)v)\, dF_V(v)$ a.s. in $\omega$, where $F_V$ is the distribution function of $V$.

Therefore, since $f_n$ is $\mathcal{G}_{n-1}$-measurable and $V_n$ is independent from $\mathcal{G}_{n-1} \vee \sigma(X_n)$, for every Borel real bounded measurable function $h$, it holds

$$\mathbb{E}(h(W_n)|\mathcal{G}_{n-1}) = \int_0^1 h(f_n v)\, dv \cdot f_n + \int_0^1 h(1 - (1 - f_n)v)\, dv \cdot (1 - f_n),$$

where we have also used $\mathbb{P}(X_n = 0|\mathcal{G}_{n-1}) = \mathbb{P}(X_n = 0|\mathcal{F}_{n-1}^X)$. The changes of variables $y = f_n v$ and $z = 1 - (1 - f_n)v$ yield

$$\mathbb{E}(h(W_n)|\mathcal{G}_{n-1}) = \int_0^{f_n} h(y)\, dy + \int_{f_n}^1 h(z)\, dz = \int_0^1 h(v)\, dv.$$

Then $W_n$ is independent of $\mathcal{G}_{n-1}$ and it is uniformly distributed in $[0,1]$. The other statements follow from the equalities

$$(2) \quad X_n = \mathbf{1}(W_n > f_n) \text{ and } V_n = \frac{W_n}{f_n}\mathbf{1}(W_n \leq f_n) + \frac{1-W_n}{1-f_n}\mathbf{1}(W_n > f_n). \qquad \square$$

Lemma 2 shows that $\mathcal{G}$ is the natural filtration of $(X,W)$ and that $(W_n : n \leq 0)$ is a sequence of independent increments for this filtration. Thus, it is direct to prove that $\mathcal{G} = \mathcal{F}^W \Leftrightarrow \mathcal{G}_0 = \mathcal{F}_0^W \Leftrightarrow \mathcal{F}_0^X \subseteq \mathcal{F}_0^W$. Therefore, $\mathcal{F}_0^X \subseteq \mathcal{F}_0^W$ is a sufficient condition for $\mathcal{G} = \mathcal{F}^W$ to be of product type, and thus, for $\mathcal{F}^X$ to be standard.

Now, the condition $\mathcal{F}_0^X \subseteq \mathcal{F}_0^W$ is not always fulfilled, even if the tail $\sigma$-field $\mathcal{F}_{-\infty}^X$ is trivial. This is one of the main points in the theory of standardness. A historical reference on this matter, that we ought to the referee, is [11], Section III, paragraph 12. In the next section we exhibit a class of processes verifying $\mathcal{F}_0^X \subseteq \mathcal{F}_0^W$.

**3. Stationary processes of summable memory decay are standard.** For $N \leq K \leq 0$, we set $X[N;K] = (X_n : n = N,\ldots,K)$ and $X(-\infty;K] = (X_n : n \leq K)$. We put $\Sigma^{(K)} = \prod_{n \leq K}\{0,1\}$, for every $K \leq 0$. A point in $\Sigma^{(K)}$ will be denoted simply by $\mathbf{x}$.

The conditional probability is written $\mathbb{P}(i|\mathbf{x}) = \mathbb{P}(X_0 = i|X(-\infty;-1] = \mathbf{x})$ for $i \in \{0,1\}$, $\mathbf{x} \in \Sigma^{(-1)}$. We assume all the cylinder sets have strictly positive measure and that $\mathbb{P}(i|\mathbf{x}) > 0$ for every $i \in \{0,1\}$, $\mathbf{x} \in \Sigma^{(-1)}$.

For $p \geq 0$, define the following quantity:

$$(3) \quad \gamma_p = 1 - \inf\left\{\frac{\mathbb{P}(i|\mathbf{x})}{\mathbb{P}(i|\mathbf{y})} : i \in \{0,1\}, \mathbf{x},\mathbf{y} \in \Sigma^{(-1)}, \mathbf{x}[-p;-1] = \mathbf{y}[-p;-1]\right\},$$

where in the case $p = 0$ there is no restriction on the variables $\mathbf{x}, \mathbf{y} \in \Sigma^{(-1)}$. The sequence $(\gamma_p : p \geq 0)$ is decreasing and $[0,1]$ valued. This process is said to have complete connections if it verifies $\lim_{p \to \infty} \gamma_p = 0$ (see [9]). Let us show that in this case $\gamma_p \in [0,1)$ for all $p \geq 0$. Simply note that if $\gamma_p < 1$ for some $p$, then $\gamma_0 < 1$, thus, $\gamma_q < 1$ for all $q$. Indeed, fix $\mathbf{v} \in \Sigma^{(-p-1)}$. Then for every $\mathbf{x}, \mathbf{y} \in \Sigma^{(-1)}$

$$\mathbb{P}(i|\mathbf{x}) \geq (1-\gamma_p)\mathbb{P}(i|\mathbf{v}\mathbf{x}[-p,-1])$$
$$\geq c =: (1-\gamma_p)\inf\{\mathbb{P}(j|\mathbf{v}z) : j \in \{0,1\}, z \in \{0,1\}^p\} > 0,$$

thus, $\frac{\mathbb{P}(i|\mathbf{x})}{\mathbb{P}(i|\mathbf{y})} \geq c$ from where we deduce $\gamma_0 \leq 1 - c$.

If the additional property $\sum_{p \geq 0} \gamma_p < \infty$ holds, the process is said to have summable memory decay. Our next result assumes a weaker condition than summable memory decay.

THEOREM 3. *Assume the process $(X_n : n \leq 0)$ has complete connections. If*

$$\sum_{\ell=0}^{\infty} \prod_{p=0}^{\ell}(1-\gamma_p) = \infty,$$

*then the filtration $\mathcal{F}^X$ is standard.*

PROOF. First, let us fix a generating r.v. $R$, that is, such that $\mathcal{F}_0^X = \sigma(R)$. We choose

(4) $$R = \sum_{n \leq 0} 3^n X_n,$$

so that, for $n \leq 0$, $\{R(\omega) - R(\omega') < 3^n\} = \{X[n;0](\omega) = X[n;0](\omega')\}$. As we pointed out, a sufficient condition ensuring $\mathcal{F}^X$ is standard is that $R$ is $\mathcal{F}_0^W$-measurable. In the sequel, for all $N \leq 0$, we will construct a function $F_N : [0,1]^{|N|+1} \to \mathbb{R}$ such that $S_N = F_N(W[N;0])$ converges in probability toward $R$, and the result will be shown.

Let us consider the sequences $(V_n : n \leq 0)$ and $(W_n : n \leq 0)$ introduced in Sections 1 and 2, so

(5) $$X_n = \mathbf{1}(W_n > \mathbb{P}(0|X(-\infty;n-1])).$$

For all $N \leq 0$, let us construct an approximation $(\widehat{X}_n^{(N)} : n \leq 0)$ of the process. Before $N$, we put (arbitrarily) $\widehat{X}_n^{(N)} = 0$ for $n < N$, and for $n \in \{N, \ldots, 0\}$, the evolution of $\widehat{X}^{(N)}$ is governed by the recurrence

(6) $$\widehat{X}_n^{(N)} = \mathbf{1}(W_n > \mathbb{P}(0|\widehat{X}^{(N)}(-\infty;n-1])).$$

We define $S_N = \sum_{n \leq 0} 3^n \widehat{X}_n^{(N)}$, then $S_N$ is a function of $W[N;0]$. To prove the theorem, it is enough to show convergence in probability of $S_N$ toward $R$. For that purpose, fix $\varepsilon > 0$ and $K$ a positive integer such that $3^{-K} < \varepsilon$. For $N$ smaller than $-K$, one has

$$\mathbb{P}(|S_N - R| > \varepsilon) \leq \mathbb{P}(|S_N - R| \geq 3^{-K}) = \mathbb{P}(\widehat{X}^{(N)}[-K;0] \neq X[-K;0]).$$

Therefore, the result will follow once we prove

(7) $$\lim_{N \to -\infty} \mathbb{P}(\widehat{X}^{(N)}[-K;0] \neq X[-K;0]) = 0.$$

The proof relies on ingredients that have been developed in [2], as well as in [5], in alternative shapes. For $i \in \{0,1\}$, set

(8) $a_0(i) = \inf\{\mathbb{P}(i|\mathbf{x}) : \mathbf{x} \in \Sigma^{(-1)}\},$

(9) $a_p(i|z) = \inf\{\mathbb{P}(i|\mathbf{x}) : \mathbf{x} \in \Sigma^{(-1)}, \mathbf{x}[-p;-1] = z\}$ for $p \geq 1$, $z \in \{0,1\}^p$.

Notice that, for all $p \geq 0$, $z \in \{0,1\}^p$ and $\mathbf{x} \in \Sigma^{(-1)}$, with $\mathbf{x}[-p;-1] = z$, it holds

(10) $\quad a_p(0|z) + a_p(1|z) \geq (1 - \gamma_p)\mathbb{P}(0|\mathbf{x}) + (1 - \gamma_p)\mathbb{P}(1|\mathbf{x}) \geq (1 - \gamma_p)$

[for $p = 0$, it simply reads $a_0(0) + a_0(1) \geq 1 - \gamma_0$].

Let $(Z_q : q \geq 0)$ be a Markov chain, taking values in $\mathbb{N}$, with initial value $Z_0 = 0$ and with transition probabilities

$$p_{i,i+1} = 1 - \gamma_i, \qquad p_{i,0} = \gamma_i, \qquad p_{i,j} = 0 \quad \text{in other cases.}$$

The hypothesis of the theorem is equivalent to the transience or null recurrence of this chain. Thus,

$$\lim_{q \to \infty} P(Z_q \leq K) = 0.$$

To prove (7), and therefore the theorem, is enough to prove the inequality

$$\mathbb{P}(\widehat{X}^{(N)}[-K;0] \neq X[-K;0]) \leq P(Z_{-N} \leq K).$$

For the rest of the proof, we follow the simplification made by the referee to our original proof. The referee introduced for $n \in \{N, \ldots, 0\}$ the random variable $L_n^{(N)} = \max\{l \in \mathbb{N} : \widehat{X}^{(N)}[n-l+1;n] = X[n-l+1;n]\}$. Notice that $\{L_0^{(N)} \leq K\} = \{\widehat{X}^{(N)}[-K;0] \neq X[-K;0]\}$.

For $n \in \{N+1, \ldots, 0\}$, it follows from the definition of $L^{(N)}$, (5) and (6) that

$$\{L_{n-1}^{(N)} = l, L_n^{(N)} = l+1\} \supseteq \{L_{n-1}^{(N)} = l, W_n < a_l(0|X[n-l;n-1])\}$$
$$\cup \{L_{n-1}^{(N)} = l, W_n > 1 - a_l(1|X[n-l;n-1])\}.$$

Thus, on the set $\{L_{n-1}^{(N)} = l\}$ we have the inequality

$$\mathbb{P}(L_n^{(N)} = l+1|\mathcal{G}_{n-1}) \geq a_l(0|X[n-l;n-1]) + a_l(1|X[n-l;n-1]) \geq 1 - \gamma_l,$$

which proves that

$$\mathbb{P}(L_n^{(N)} = L_{n-1}^{(N)} + 1|\mathcal{G}_{n-1}) \geq 1 - \gamma_{L_{n-1}^{(N)}}.$$

Now, let us prove by induction on $n \in \{N, \ldots, 0\}$ that $L_n^{(N)} \geq Z_{n-N}$ in law, namely,

(11) $\qquad \mathbb{P}(L_n^{(N)} > M) \geq \mathbb{P}(Z_{n-N} > M) \qquad \text{for all } M \in \mathbb{N}.$

For $n = N$, this is obvious because $Z_0 = 0$. Assuming the inequality holds for a given $n \leq -1$, we get

$$\mathbb{P}(L_{n+1}^{(N)} > M) = \mathbb{P}(L_n^{(N)} \geq M, L_{n+1}^{(N)} = L_n^{(N)} + 1)$$

$$\geq \mathbb{E}(\mathbf{1}(L_n^{(N)} \geq M)(1 - \gamma_{L_n^{(N)}}))$$

$$\geq \mathbb{E}(\mathbf{1}(Z_{n-N} \geq M)(1 - \gamma_{Z_{n-N}}))$$

$$= \mathbb{P}(Z_{n-N} \geq M, Z_{n-N+1} = Z_{n-N} + 1)$$

$$= \mathbb{P}(Z_{n-N+1} > M).$$

Here we have used that $L_n^{(N)} \geq Z_{n-N}$, in law, and that the function $l \to \mathbf{1}(l \geq M)(1 - \gamma_l)$ is increasing. The theorem is finally obtained by taking $n = 0$ in (11). $\square$

REMARK 4. We notice that if $\gamma_p = 0$ for some $p \geq 1$, the process $((X_{n-p+1}, \ldots, X_n) : n \leq 0)$ is a Markov chain and Theorem 3 is well known (see [12]). When $p = 0$, the result is trivial because $(X_n : n \leq 0)$ are independent.

**4. A product type filtration assuming standardness.** In this section we assume $\mathcal{F}^X$ is standard. As stated, we will construct a diffusive product type extension of $\mathcal{F}^X$. We consider the sequences $(V_n : n \leq 0)$ and $(W_n : n \leq 0)$ introduced in Sections 1 and 2, and the filtration $\mathcal{G} = (\mathcal{G}_n : n \leq 0)$ defined by $\mathcal{G}_n = \sigma(X_m, V_m : m \leq n)$. For a notational purpose, if $Z$ and $Z'$ are random elements, we denote by $\mathcal{L}(Z)$ the probability distribution of $Z$ and by $\mathcal{L}(Z|Z' = z')$ its conditional law with respect to the event $\{Z' = z'\}$.

Let $\rho_0$ be a metric in $\Sigma$, consider the following sequence $(\rho_{|n|} : n \leq 0)$ defined recursively, for $n \leq -1$ and $\mathbf{x}, \mathbf{y} \in \Sigma$, by

$$\rho_{|n|}(\mathbf{x}, \mathbf{y})$$
$$(12) \quad = \inf\{\mathbb{E}_\Lambda(\rho_{|n|-1}(\mathbf{x}(-\infty; n]\xi 0^{|n|-1},$$
$$\mathbf{y}(-\infty; n]\eta 0^{|n|-1})) : \Lambda \in \mathcal{J}(\mathbf{x}(-\infty; n], \mathbf{y}(-\infty; n])\},$$

where, for every $\mathbf{z}, \mathbf{w} \in \Sigma$, $\mathcal{J}(\mathbf{z}, \mathbf{w})$ is the set of couplings of $\xi$ and $\eta$ whose marginals satisfy $\mathcal{L}(\xi) = \mathcal{L}(X_{n+1}|X(-\infty; n] = \mathbf{z})$ and $\mathcal{L}(\eta) = \mathcal{L}(X_{n+1}|X(-\infty; n] = \mathbf{w})$. We have put $0^{|n|-1} = \underbrace{0 \ldots 0}_{|n|-1 \text{ times}}$, but instead of $0^{|n|-1}$, any other fixed choice can also be taken.

If $\mathcal{F}^X$ is standard, it satisfies Vershik criterion (see [15, 16]): for all initial metric $\rho_0$,

$$(13) \quad \lim_{p \to \infty} \alpha_p(\rho_0) = 0 \quad \text{where } \alpha_p(\rho_0) = \int_{\Sigma \times \Sigma} \rho_p(\mathbf{x}, \mathbf{y}) \, d\mathbb{P}(\mathbf{x}) \, d\mathbb{P}(\mathbf{y})$$

$$\text{for } p \geq 0.$$

From the cosiness property introduced in [14] (see also [6, 7, 10]), it suffices to verify (13) for the following well-defined metric $\rho_0(\mathbf{x}, \mathbf{y}) = |R(\mathbf{x}) - R(\mathbf{y})|$, for a generating function $R$. We point out that, in the case of stationary processes, this property will also follow from our construction. We fix $R$ as in (4), and our construction will depend on this arbitrary choice.

From its definition, $\rho_{|n|}(\mathbf{x},\mathbf{y})$ does not depend on $(\mathbf{x}[n+1;0], \mathbf{y}[n+1;0])$, so, since the process is stationary, we get $\alpha_{|n|}(\rho_0) = \int_{\Sigma \times \Sigma} \widetilde{\rho}_{|n|}(\mathbf{x},\mathbf{y})\,d\mathbb{P}(\mathbf{x})\,d\mathbb{P}(\mathbf{y})$, where we set $\widetilde{\rho}_{|n|}(\mathbf{x},\mathbf{y}) = \rho_{|n|}(\mathbf{x}0^{|n|}, \mathbf{y}0^{|n|})$.

For $\mathbf{x},\mathbf{y} \in \Sigma^{(-1)}$, consider

$$\lambda_m(\mathbf{x},\mathbf{y}) = \mathrm{sign}(\widetilde{\rho}_{|m|-1}(\mathbf{x}0,\mathbf{y}0)$$
$$+ \widetilde{\rho}_{|m|-1}(\mathbf{x}1,\mathbf{y}1) - \widetilde{\rho}_{|m|-1}(\mathbf{x}0,\mathbf{y}1) - \widetilde{\rho}_{|m|-1}(\mathbf{x}1,\mathbf{y}0)).$$

A direct computation shows that the following coupling minimizes the expectation $\mathbb{E}_\Lambda(\widetilde{\rho}_{|m|-1}(\mathbf{x}\xi, \mathbf{y}\eta))$:

| $\xi \setminus \eta$ | 0 | 1 | |
|---|---|---|---|
| 0 | $\mathbb{P}(0\|\mathbf{x}) \wedge \mathbb{P}(0\|\mathbf{y})$ | $(\mathbb{P}(0\|\mathbf{x}) - \mathbb{P}(0\|\mathbf{y}))^+$ | if $\lambda_m(\mathbf{x},\mathbf{y}) = -1$ |
| 1 | $(\mathbb{P}(1\|\mathbf{x}) - \mathbb{P}(1\|\mathbf{y}))^+$ | $\mathbb{P}(1\|\mathbf{x}) \wedge \mathbb{P}(1\|\mathbf{y})$ | |

and

| $\xi \setminus \eta$ | 0 | 1 | |
|---|---|---|---|
| 0 | $(\mathbb{P}(0\|\mathbf{x}) - \mathbb{P}(1\|\mathbf{y}))^+$ | $\mathbb{P}(0\|\mathbf{x}) \wedge \mathbb{P}(1\|\mathbf{y})$ | if $\lambda_m(\mathbf{x},\mathbf{y}) = 1$ |
| 1 | $\mathbb{P}(1\|\mathbf{x}) \wedge \mathbb{P}(0\|\mathbf{y})$ | $(\mathbb{P}(1\|\mathbf{x}) - \mathbb{P}(0\|\mathbf{y}))^+$ | |

(see [4], Lemma 5.2, for a similar construction). This coupling is denoted by $\Lambda_m(\cdot,\cdot|\mathbf{x},\mathbf{y}) \in \mathcal{J}(\mathbf{x},\mathbf{y})$.

With this notation, we can write $\rho_{|n|}$ in terms of $\rho_{|n|-1}$ by

$$(14) \quad \rho_{|n|}(\mathbf{x},\mathbf{y}) = \mathbb{E}_{\Lambda_n(\cdot,\cdot|\mathbf{x},\mathbf{y})}(\rho_{|n|-1}(\mathbf{x}(-\infty;n]\xi 0^{|n|-1}, \mathbf{y}(-\infty;n]\eta 0^{|n|-1})).$$

For each fixed $N \leq 0$ and a point $\widehat{\mathbf{x}}^{(N)} \in \Sigma$, we construct an approximation $\widehat{X}^{(N)}[N;0]$ of $X[N;0]$ and a sequence $U^{(N)}[N;0]$ of uniform i.i.d. r.v.'s, defined recursively and such that $\widehat{X}^{(N)}[N;0]$ is measurable with respect to $\sigma(U^{(N)}[N;0])$. This is done inductively starting with $\widehat{X}^{(N)}(-\infty;N-1] = \widehat{\mathbf{x}}^{(N)}(-\infty;N-1]$.

DEFINITION 5. Consider $m \in \{N-1,\ldots,-1\}$ and define

$$(15) \quad U_{m+1}^{(N)} = \begin{cases} W_{m+1}, & \text{on } \lambda_m(X(-\infty;m], \widehat{X}^{(N)}(-\infty;m]) = -1, \\ 1 - W_{m+1}, & \text{on } \lambda_m(X(-\infty;m], \widehat{X}^{(N)}(-\infty;m]) = 1, \end{cases}$$

and

$$(16) \quad \widehat{X}_{m+1}^{(N)} = \mathbf{1}(U_{m+1}^{(N)} > \mathbb{P}(0|\widehat{X}^{(N)}(-\infty;m])).$$

In the sequel we specify the structure of the sequence and explain how to recover $X$ from $U^{(N)}$. We also study the joint law of $X$ and $\widehat{X}^{(N)}$.

LEMMA 6. *$U^{(N)}[N;0]$ is a sequence of i.i.d. r.v.'s uniformly distributed on $[0,1]$. For all $m \in \{N,\ldots,0\}$, $U_m^{(N)}$ is independent of $\mathcal{G}_{m-1}$. Moreover, $\mathcal{G}_{m-1} \vee \sigma(U_m^{(N)}) = \mathcal{G}_m$.*

PROOF. Let $m \in \{N,\ldots,0\}$. The law of $U_m^{(N)}$ given $\mathcal{G}_{m-1}$ is the same as the law of $W_m$ given $\mathcal{G}_{m-1}$. Then, the uniform distribution of $U_m^{(N)}$ on $[0,1]$ and the independence between $U_m^{(N)}$ and $\mathcal{G}_{m-1}$ readily follow.

To conclude, let us express explicitly $X_m$ in terms of $X(-\infty;m-1]$, $\widehat{X}(-\infty;m-1]$ and $U_m^{(N)}$. From (1) and (15), we get the following:

- if $\lambda_{m-1}(X(-\infty;m-1],\widehat{X}^{(N)}(-\infty;m-1]) = -1$, then $X_m = \mathbf{1}(U_m^{(N)} > \mathbb{P}(0|X(-\infty;m-1]))$,
- if $\lambda_{m-1}(X(-\infty;m-1],\widehat{X}^{(N)}(-\infty;m-1]) = 1$, then $X_m = \mathbf{1}(1 - U_m^{(N)} > \mathbb{P}(0|X(-\infty;m-1]))$,

where $\widehat{X}^{(N)}(-\infty;m-1]$ is itself a function of $X(-\infty;m-1]$, $U^{(N)}[N;m-1]$ and $\widehat{\mathbf{x}}^{(N)}(-\infty,N-1]$. □

We observe that $\mathbb{P}(\widehat{X}_m^{(N)} = 0) = \mathbb{P}(0|\widehat{X}^{(N)}(-\infty;m-1])$. Finer relations are given in Lemma 7 below.

Let us write how to recover the whole sequence $X[N;0]$ from $U^{(N)}[N;0]$ and the past. We define a function $G: \{1,-1\} \times [0,1] \times \Sigma \to \{0,1\}$ by

$$G(\lambda, u, \mathbf{x}) = \begin{cases} \mathbf{1}(u > \mathbb{P}(0|\mathbf{x})), & \text{if } \lambda = -1, \\ \mathbf{1}(1-u > \mathbb{P}(0|\mathbf{x})), & \text{if } \lambda = 1. \end{cases}$$

We get $X_m = G(\lambda_{m-1}(X(-\infty;m-1],\widehat{X}^{(N)}(-\infty;m-1]), U_m^{(N)}, X(-\infty;m-1])$. Iterating this procedure, we can define functions $G_N$, such that

(17) $$X[N;0] = G_N(U^{(N)}[N;0], X(-\infty;N-1]).$$

We notice that $\widehat{X}^{(N)}[N;0]$ is a similar function of $U^{(N)}[N;0]$ and $\widehat{\mathbf{x}}^{(N)}(-\infty,N-1]$ (but simpler, in the sense that it does not use $\lambda$, or, equivalently, this corresponds to $\lambda_m(\widehat{X}^{(N)}(-\infty;m],\widehat{X}^{(N)}(-\infty;m]) = -1$).

LEMMA 7. *For any sequence $\mathbf{a} \in \Sigma$,*

$$\mathbb{P}(\widehat{X}^{(N)}[N;0] = \mathbf{a}[N;0])$$
$$= \mathbb{P}(X[N;0] = \mathbf{a}[N;0]|X(-\infty;N-1] = \widehat{\mathbf{x}}^{(N)}(-\infty;N-1]).$$

*For all $m \in \{N,\ldots,0\}$, and all $a,b \in \{0,1\}$,*

(18)
$$\mathbb{P}(X_m = a, \widehat{X}_m^{(N)} = b|\mathcal{G}_{m-1})$$
$$= \Lambda_{m-1}(a,b|X(-\infty;m-1],\widehat{X}^{(N)}(-\infty;m-1]).$$

PROOF. Let us write the joint law $\mathcal{L}(X_m, \widehat{X}_m^{(N)}|\mathcal{G}_{m-1})$. Since $\lambda_{m-1}(X(-\infty; m-1], \widehat{X}^{(N)}(-\infty; m-1])$ is $\mathcal{G}_{m-1}$-measurable, we can treat the cases according to the values of this variable. We only check one case, $(a,b) = (0,0)$ and $\lambda_{m-1}(X(-\infty; m-1], \widehat{X}^{(N)}(-\infty; m-1]) = -1$. One has

$$\mathbb{P}(X_m = 0, \widehat{X}_m^{(N)} = 0|\mathcal{G}_{m-1})$$
$$= \mathbb{P}(W_m \leq \mathbb{P}(0|\widehat{X}^{(N)}(-\infty; m-1])|X_m = 0, \mathcal{G}_{m-1})\mathbb{P}(X_m = 0|\mathcal{G}_{m-1})$$
$$= \mathbb{P}(\mathbb{P}(0|X(-\infty; m-1])V_m \leq \mathbb{P}(0|\widehat{X}^{(N)}(-\infty; m-1])|X_m = 0, \mathcal{G}_{m-1})$$
$$\quad \times \mathbb{P}(0|X(-\infty; m-1])$$
$$= \mathbb{P}(0|X(-\infty; m-1]) \wedge \mathbb{P}(0|\widehat{X}^{(N)}(-\infty; m-1]),$$

where the last line follows since $V_m$ is a uniform random variable independent of $\mathcal{G}_{m-1} \vee \sigma(X_m)$. □

We define $\widehat{R}^{(N)} = R(\widehat{X}^{(N)}(-\infty; 0])$. Therefore, $\widehat{R}^{(N)}$ is generated by the sequence $U^{(N)}[N; 0]$ and it is independent of $X(-\infty; N-1]$.

LEMMA 8. *The following equality holds:* $\mathbb{E}(|R - \widehat{R}^{(N)}|) = \int_\Sigma \rho_{|N|+1}(\mathbf{x}, \widehat{\mathbf{x}}^{(N)}) d\mathbb{P}(\mathbf{x})$.

PROOF. We must show $\mathbb{E}(\rho_0(X, \widehat{X}^{(N)})) = \int_\Sigma \rho_{|N|+1}(\mathbf{x}, \widehat{\mathbf{x}}^{(N)}) d\mathbb{P}(\mathbf{x})$. Notice that $\rho_{|N|+1}$ does not depend on coordinates $\{N, \ldots, 0\}$, so

$$\int_\Sigma \rho_{|N|+1}(\mathbf{x}, \widehat{\mathbf{x}}^{(N)}) d\mathbb{P}(\mathbf{x})$$
$$= \mathbb{E}(\rho_{|N|+1}(X, \widehat{\mathbf{x}}^{(N)}))$$
$$= \mathbb{E}(\rho_{|N|+1}(X(-\infty; N-1]0^{|N|+1}, \widehat{X}^{(N)}(-\infty; N-1]0^{|N|+1})).$$

Recall (14), that in our case reads, for $m \leq -1$,

$$\rho_{|m|}(X(-\infty; m]0^{|m|}, \widehat{X}^{(N)}(-\infty; m]0^{|m|})$$
$$= \mathbb{E}_{\Lambda_m(\cdot, \cdot|X(-\infty; m], \widehat{X}^{(N)}(-\infty; m])}$$
$$\quad \times (\rho_{|m|-1}(X(-\infty; m]\xi 0^{|m|-1}, \widehat{X}^{(N)}(-\infty; m]\eta 0^{|m|-1})).$$

Then, Lemma 7 shows that, for any measurable function $h$, it holds:

$$\mathbb{E}(\mathbb{E}_{\Lambda_m(\cdot, \cdot|X(-\infty; m], \widehat{X}^{(N)}(-\infty; m])}(h(X(-\infty; m]\xi, \widehat{X}^{(N)}(-\infty; m]\eta)))$$
$$= \mathbb{E}(h(X(-\infty; m+1], \widehat{X}^{(N)}(-\infty; m+1])).$$

Hence,
$$\mathbb{E}(\rho_{|m|}(X(-\infty;m]0^{|m|}, \widehat{X}^{(N)}(-\infty;m]0^{|m|}))$$
$$= \mathbb{E}(\rho_{|m|-1}(X(-\infty;m+1]0^{|m|-1}, \widehat{X}^{(N)}(-\infty;m+1]0^{|m|-1})).$$

The argument holds for all $m \in \{N-1, \ldots, -1\}$ and the lemma is proved. □

$R$ is determined from the whole past up to $N-1$ and the i.i.d. r.v.'s $U^{(N)}[N;0]$. In fact, from (17), $R(X(-\infty;0]) = R(X(-\infty;N-1]G_N(U^{(N)}[N;0], X(-\infty;N-1]))$.

The following result is a direct consequence of the martingale theorem, and we skip a detailed proof.

LEMMA 9. *Let $N \leq 0$, $\delta > 0$, $Z[N;0]$ be a sequence of uniform i.i.d. r.v. independent of $X(-\infty;N-1]$ and $H$ a measurable function such that*
$$X[N;0] = H(Z[N;0], X(-\infty;N-1]).$$
*Then, there exists an integer $K = K(N, \delta, H) < N$ and a function $\Phi : [0,1]^{|N|+1} \times \{0,1\}^{N-K} \to \mathbb{R}$, which depends on $N, \delta, H$, that verify*
$$\mathbb{P}(|\Phi(Z[N;0], X[K;N-1]) - R| > \delta) < \delta.$$

One of the tools we need is given by the following construction. Let us take $\delta > 0$ and consider $N = N(\delta) \leq 0$ such that $\alpha_{|N|+1}(\rho_0) < \delta$. By Fubini's theorem, we can choose a sequence $\widehat{\mathbf{x}}^{(N)} \in \Sigma$ verifying the following property:

$$(19) \qquad \int_\Sigma \rho_{|N|+1}(\mathbf{x}, \widehat{\mathbf{x}}^{(N)}) \, d\mathbb{P}(\mathbf{x}) < \delta.$$

The choice of such $\widehat{\mathbf{x}}^{(N)}$ for each relevant $N$ is arbitrary and will influence our construction. From Lemma 8, we obtain that, for such $N$ and $\widehat{\mathbf{x}}^{(N)}$, the next bound holds:
$$\mathbb{E}(|R - \widehat{R}^{(N)}|) \leq \delta.$$

Now we construct a sequence $(U_n : n \leq 0)$ of uniform i.i.d. r.v. that will give us a product type filtration such that $\mathcal{F}^X$ is immersed on. Fix a positive sequence $(\delta_j : j \geq 0)$ decreasing to 0.

- Initially, at step 0, we choose $N_0$ and $\widehat{\mathbf{x}}^{(N_0)} \in \Sigma$ such that $\alpha_{|N_0|+1}(\rho_0) < \delta_0$ and
$$\int \rho_{|N_0|+1}(\mathbf{x}, \widehat{\mathbf{x}}^{(N_0)}) \, d\mathbb{P}(\mathbf{x}) < \delta_0.$$
We construct $U^{(N_0)}[N_0;0]$ and $\widehat{X}^{(N_0)}[N_0;0]$ following Definition 5. We put $M_0 = 1$, $M_1 = N_0$ and $H_0 = G_{N_0}$, so that $X[M_1;0] = H_0(U^{(N_0)}[M_1;0], X(-\infty;M_1-1])$, see (17). In particular, we have that $\mathbb{E}(|R - \widehat{R}^{(N_0)}|) \leq \delta_0$. We finally put $U[N_0;0] = U^{(N_0)}[N_0;0]$.

- Assume at step $j-1$ we have constructed a sequence $U[M_j;0]$ and a function $H_{j-1}$ such that

(20) $$X[M_j;0] = H_{j-1}(U[M_j;0], X(-\infty; M_j - 1]).$$

We obtain $K_j < M_j$ and $\Phi_j$ by applying Lemma 9 with $N = M_j$, $\delta = \delta_j/2$, $Z[M_j;0] = U[M_j;0]$ and $H = H_{j-1}$. We choose $N_j$ and $\widehat{\mathbf{x}}^{(N_j)}$ such that

(21) $$\alpha_{|N_j|+1}(\rho_0) < 3^{K_j - M_j + 1} \cdot \delta_j/2 \quad \text{and}$$

$$\int \rho_{|N_j|+1}(\mathbf{x}, \widehat{\mathbf{x}}^{(N_j)}) \, d\mathbb{P}(\mathbf{x}) < 3^{K_j - M_j + 1} \cdot \delta_j/2.$$

We set $M_{j+1} = M_j + N_j - 1$.

- Applying the construction on the shifted process $(X_{n+M_j-1} : n \leq 0)$ and using stationarity, we construct a sequence $U[M_{j+1}; M_j - 1]$ of uniform i.i.d. r.v., which is independent of $U[M_j;0]$, such that

(22) $$X[M_{j+1}; M_j - 1] = G_{N_j}(U[M_{j+1}; M_j - 1], X(-\infty; M_{j+1} - 1]).$$

From (20) and (22), we can define a function $H_j$ in terms of $G_{N_j}$ and $H_{j-1}$ such that $X[M_{j+1};0] = H_j(U[M_{j+1};0], X(-\infty; M_{j+1} - 1])$.

A repeated use of Lemma 6 in the construction of the blocks $U[M_{j+1}; M_j - 1]$ gives that $(U_n : n \leq 0)$ is a sequence of i.i.d. r.v.'s uniformly distributed in $[0,1]$, so $\mathcal{F}^U$ is a diffusive product type filtration.

THEOREM 10. *If $\mathcal{F}^X$ is standard, then $\mathcal{F}^X$ is immersed in the diffusive product type filtration $\mathcal{F}^U$.*

PROOF. It is enough to construct a function $S$ such that $R(X(-\infty;0]) = S(U(-\infty;0])$. For $j \geq 1$, set $S_j(w) = \Phi_j(U[M_j;0](w), \widehat{X}[K_j; M_j - 1](w))$, where $\widehat{X} = \widehat{X}^{(M_{j+1})}$ is the process generated in Definition 5 starting from $\widehat{\mathbf{x}}^{(N_j)}$. This means $\widehat{X}(-\infty; M_{j+1} - 1] = \widehat{\mathbf{x}}^{(N_j)}(-\infty; N_j - 1]$, where we identify points in $\Sigma^{(M_{j+1}-1)}$ and $\Sigma^{(N_j-1)}$. Therefore, $S_j$ is a function of $U[M_{j+1};0]$ because $\widehat{X}[K_j; M_j - 1]$ is a function of $U[M_{j+1}; M_j - 1]$. It remains to prove that $S_j$ converges in probability to $R$.

Notice that $X[K_j; M_j - 1] = \widehat{X}[K_j; M_j - 1]$ implies $S_j = \Phi_j(U[M_j;0], X[K_j; M_j - 1])$. Then

$$\mathbb{P}(S_j \neq \Phi_j(U[M_j;0], X[K_j; M_j - 1])) \leq P(X[K_j; M_j - 1] \neq \widehat{X}[K_j; M_j - 1]).$$

Recall that $|R(\mathbf{x}) - R(\mathbf{y})| < 3^{-k}$ implies $\mathbf{x}[-k;0] = \mathbf{y}[-k;0]$, then we get

$$\mathbb{P}(X[K_j; M_j - 1] \neq \widehat{X}[K_j; M_j - 1])$$
$$\leq \mathbb{P}(|R(X(-\infty; M_j - 1]) - R(\widehat{X}(-\infty; M_j - 1])| \geq 3^{-(M_j - 1 - K_j)})$$
$$\leq 3^{M_j - 1 - K_j} \mathbb{E}(|R(X(-\infty; M_j - 1]) - R(\widehat{X}(-\infty; M_j - 1])|),$$

where we have identified $\Sigma$ and $\Sigma^{(M_j-1)}$. By applying Lemma 8 to the shifted process and in view of the choice of $N_j$ in (21), we find

$$\mathbb{E}(|R(X(-\infty; M_j - 1]) - R(\widehat{X}(-\infty; M_j - 1])|) \leq 3^{K_j - M_j + 1}\delta_j/2.$$

We have proven $\mathbb{P}(S_j \neq \Phi_j(U[M_j; 0], X[K_j; M_j - 1])) \leq \delta_j/2$. On the other hand, the choice of $K_j$ done in Lemma 9 guarantees that $\mathbb{P}(|\Phi_j(U[M_j; 0], X[K_j; M_j - 1]) - R(X(-\infty, 0])| > \delta_j/2) \leq \delta_j/2$. Therefore,

$$\begin{aligned}
\mathbb{P}(|S_j &- R(X(-\infty, 0])| > \delta_j) \\
&\leq \mathbb{P}(S_j \neq \Phi_j(U[M_j; 0], X[K_j; M_j - 1])) \\
&\quad + \mathbb{P}(|\Phi_j(U[M_j; 0], X[K_j; M_j - 1]) - R(X(-\infty, 0])| > \delta_j/2) \leq \delta_j,
\end{aligned}$$

then the convergence in probability follows. □

**Acknowledgments.** The authors thank M. Émery for fruitful discussions on standardness during his visit to the Workshop "Potential, Probability and Filtrations" at the CMM. The authors are indebted to an anonymous referee, whose suggestions allowed us to correct some mistakes in the original version which helped us improve the whole presentation of the paper, in particular, the proof of Theorem 3, and because he/she called our attention to reference [11].

X. Bressaud
Institut Mathématiques de Luminy
Université de la Mediterranée
Marseille
France
E-mail: bressaud@iml.univ-mrs.fr

A. Maass
S. Martinez
J. San Martin
CMM-DIM
Universidad de Chile
Casilla 170-3, Correo 3
Santiago
Chile
E-mail: amaass@dim.uchile.cl
         smartine@dim.uchile.cl
         jsanmart@dim.uchile.cl